\documentclass[12pt]{amsart}
\usepackage[all, cmtip]{xy}

\usepackage{times,amsfonts,amsmath,amstext,amsbsy,amssymb, amsopn,amsthm,upref,eucal, mathrsfs}

\newcommand \la {\lambda}

\newcommand \Conf {{\mathrm {Conf}}}

\newcommand \Leb {{\mathrm {Leb}}}

\newcommand \Prob {{\mathbb P}}

\newtheorem{theorem}{Theorem}[section]

\newtheorem{lemma}[theorem]{Lemma}
\newtheorem{assumption}{Assumption}
\newtheorem{corollary}[theorem]{Corollary}
\newtheorem{proposition}[theorem]{Proposition}

\begin{document}
\title[ Conditional measures of determinantal point processes]{Conditional measures of determinantal point processes}

\author{Alexander I. Bufetov}
\address{Aix-Marseille Universit{\'e}, CNRS, Centrale Marseille, I2M, UMR 7373}\address{ The Steklov Institute of Mathematics, Moscow}\address{
The Institute for Information Transmission Problems, Moscow}\address{
National Research University Higher School of Economics, Moscow}
\begin{abstract}
For a class of one-dimensional
determinantal point processes including those  induced by orthogonal projections with integrable kernels satisfying a growth condition, it is proved that their conditional  measures, 
 with respect to the configuration in the complement of a compact interval, 
are orthogonal polynomial ensembles with explicitly found weights. Examples include the sine-process and the  process with the Bessel kernel.
The  argument uses the quasi-invariance, established in \cite{Buf-gibbs}, of our point processes under the group of piecewise isometries of ${\mathbb R}$.
\end{abstract} 
\maketitle
\section{Formulation of the main result.}
\subsection{Conditional measures.}

Let $E$ be a locally compact complete metric space, let $\Conf(E)$ be the space of configurations on $E$. Given a configuration $X \in \textrm{Conf}(E)$ and a subset $C \subset E$, we let $X|_C$ stand for the restriction of $X$ onto the subset $C$.

A  point process on $E$ is a Borel probability measure on $\Conf(E)$. 
For such a measure $\Prob$, the measure $\mathbb{P}(\cdot | X; C)$ on $\Conf(E\setminus C)$ is defined as the conditional measure of $\mathbb{P}$ with respect to the condition that the restriction of our random configuration onto $C$ coincides with $X|_C$. More formally, consider the surjective restriction mapping $X \to X|_C$ from
$\textrm{Conf}(E)$ to $\textrm{Conf}(C)$. Fibres of this mapping can be identified with $\mathrm{Conf}(E\backslash C)$, and conditional measures, in the sense of Rohlin \cite{Rohmeas},  are precisely the measures $\mathbb{P}(\cdot | X; C)$.
If the point process $\Prob$ admits correlation measures of order up to $l$, then, given distinct points $q_1, \dots, q_l\in E$, we let $\Prob^{q_1, \dots, q_l}$ stand for the $l$-th reduced Palm measure of $\Prob$ conditioned at points $q_1, \dots, q_l$ (here and below we follow the conventions of \cite{Buf-gibbs} in working with Palm measures). 

The main results of this note can informally be summarized as follows. 
If the measure $\mathbb{P}(\cdot | X; C)$ is supported on the subspace of configurations 
with precisely $l$ particles and  the reduced Palm measures, 
conditioned at different $l$-tuples of points, are equivalent, then, under certain additional assumptions (see Proposition \ref{palm-cond-char} below),  
the conditional measure $\mathbb{P}(\cdot | X; C)$ has the form
$$
Z^{-1}(q_1, \dots, q_l)\frac{d\Prob^{p_1, \dots, p_l}}{d\Prob^{q_1, \dots, q_l}}
\left( X|_C\right) d{\rho}_l(p_1, \dots, p_l),
$$
where $q_1, \dots, q_l$ is almost any fixed $l$-tuple, $\rho_l$ is the $l$-th correlation measure of $\Prob$  and $Z(q_1, \dots, q_l)$ is the normalization constant.
In particular, for one-dimensional determinantal processes induced by projections with integrable kernels satisfying a growth condition and $C$  the complement of a compact interval, 
it is  proved that $\mathbb{P}(\cdot | X; C)$ is an orthogonal polynomial ensemble 
with the weight found explicitly. We proceed to precise formulations.

Given a compact subset $B \subset E$ and a configuration $X \in \text{Conf}(E),$ let $ \#_B(X)$ stand for the number of particles of $X$ lying in $B$. Given a Borel subset $C \subset E,$ we let $\mathcal{F}_C$ be the $\sigma$-algebra generated by all random variables of the form $\#_B, B\subset C.$ Write $\mathcal{F}_C^{\mathbb{P}}$ for the $\mathbb{P}$-completion of $\mathcal{F}_C$.

	{\bf Definition}  (Ghosh and Peres \cite{Ghosh-sine}, \cite{Ghosh-rigid}). A point process $\mathbb{P}$ on $E$ is called {\it rigid} if for any compact  subset $B \subset E$ the function $\#_B$ is  $\mathcal{F}_{E\backslash B}^{\mathbb{P}}$-measurable.

For a subset $C\subset E$ and a natural number $l$, we write $\Conf_l(C)$ for the space of $l$-particle configurations on $C$; in other words, the space of all subsets of $C$ of cardinality $l$.
Rigidity implies that for  any precompact set $B\subset E$ and $\Prob$-almost any $X$ the conditional measure 
$\Prob(\cdot |X; E\setminus B)$ is supported on the subset $\Conf_l(B)$, where $l=\#_B(X)$.

 Let  $U\subset {\mathbb R}$ be an open set endowed with the Lebesgue measure $\Leb$.
Let $\Pi(x,y)$, $x,y\in U$, be a kernel smooth in the totality of variables. Assume that the kernel $\Pi$ induces an operator of orthogonal projection acting in $L_2(U, \Leb)$;
slightly abusing notation, we keep the same symbol $\Pi$ for this operator.
Let $L$ be the range of $\Pi$. By the Macchi-Soshnikov Theorem, the operator $\Pi$ induces a
determinantal measure $\Prob_{\Pi}$ on $\Conf(U)$.

For the sine-process, rigidity is due to Ghosh \cite{Ghosh-sine};
for the determinantal point processes with the Airy and the Bessel kernel, rigidity has been established in \cite{Buf-rigid}.

For $p\in U$, set $L(p)=\{\varphi\in L: \varphi(p)=0\}$
 and let $\Pi^p$ be the operator of orthogonal projection
onto $L(p)$. By the Shirai-Takahashi Theorem \cite{ShirTaka1}, the determinantal measure $\Prob_{\Pi^p}$ induced by the operator $\Pi^p$ is the reduced Palm measure of $\Prob_{\Pi}$ at the point $p$:
$
\Prob_{\Pi^p}=\Prob_{\Pi}^p.
$

\begin{assumption}\label{rn-xpxq}
Let $p\in U$. If $\varphi\in L$ is such that $\varphi(p)=0$, then
$
\displaystyle \frac{\varphi(x)}{x-p}\in L.
$
\end{assumption}
Proposition 3.3 in \cite{Buf-gibbs} shows  that Assumption \ref{rn-xpxq} holds, in particular, for
kernels $\Pi$ having integrable form
$
\Pi(x,y)=\displaystyle \frac{A(x)B(y) - B(x)A(y)}{x-y}.
$
\subsection{The trace-class case}
In the first theorem, we will make a restrictive
\begin{assumption}\label{as-trace}
We have 
$
\displaystyle \int\limits_U \frac{\Pi(x,x)dx}{1+|x|}<+\infty.
$
\end{assumption}
  The Bessel kernel satisfies Assumption \ref{as-trace}.
Under Assumption \ref{as-trace}, the operators $(|x|+1)^{-1}\Pi$  and
$
(x+i)^{-1}\Pi
$
belong to the trace class, and   for $p,q\in U$, the multiplicative functional 
\begin{equation}\label{mult-fun}
\Psi_{p,q}^{\Pi}(X)= \prod\limits_{x\in X}
\left(\frac{x-p}{x-q}\right)^2
\end{equation}
exists and belongs to  $L_1( \Conf(U), \Prob_{\Pi^q})$.
By Corollary 4.12 in \cite{Buf-gibbs}, 
we have the $\Prob_{\Pi}^q$-almost sure equality 
$
\displaystyle \frac{d\Prob_{\Pi}^p}{d\Prob_{\Pi}^q}=Z_{p,q}^{-1}\Psi^{\Pi}_{pq},
$
where $Z_{p,q}$ is the normalization constant. Since, for  $p,q,r\in U$, we have
$
\displaystyle \frac{d\Prob_{\Pi}^p}{d\Prob_{\Pi}^q}=\frac{d\Prob_{\Pi}^p}{d\Prob_{\Pi}^r}\frac{d\Prob_{\Pi}^r}{d\Prob_{\Pi}^q},
$ 
there exists a positive function $\rho^{\Pi}: U\to {\mathbb R}$ such that $
\displaystyle \frac{\Pi(p,p)}{\Pi(q,q)}\displaystyle \frac{d\Prob_{\Pi}^p}{d\Prob_{\Pi}^q}=\displaystyle \frac{\rho^{\Pi}(p)}{\rho^{\Pi}(q)}\Psi^{\Pi}_{pq},
$ or, equivalently, that
\begin{equation}\label{palm-rho-pi}
\int \limits_{\Conf(U)} {\Psi}^{\Pi}_{p,q}(X)d\Prob_{\Pi^q}(X)=\frac{ \rho^{\Pi}(q)}{\rho^{\Pi}(p)}\frac{\Pi(p,p)}{\Pi(q,q)}.
\end{equation}

If $\Pi$ is the Christoffel-Darboux kernel of a family of orthogonal polynomials and $\Prob_{\Pi}$ the corresponding  orthogonal polynomial ensemble, then $\rho^{\Pi}$ is the weight. The function $\rho^{\Pi}$ is defined up to a multiplicative constant. 
\begin{theorem}\label{gibbs-trace}
Let $U\subset {\mathbb R}$ be an open set. 
Let $\Pi(x,y)$, $x,y\in U$,  be a  smooth kernel that 
 induces an operator of orthogonal projection acting in $L_2(U, \Leb)$, satisfying Assumptions \ref{rn-xpxq}, \ref{as-trace} and such that the determinantal point process $\Prob_{\Pi}$ is rigid. Let $I\subset U$ be a compact interval.
Then

1.  For almost any $2l$ distinct points $p_1, \dots, p_l, q_1, \dots, q_l\in U$, we have the $\Prob^{q_1, \dots, q_l}$-almost sure
equality
$$
\frac {d\Prob^{p_1, \dots, p_l}}{d\Prob^{q_1, \dots, q_l}}(X)=
\frac {\det \Pi(q_i, q_j)|_{i,j=1, \dots, l}}{\det \Pi(p_i, p_j)|_{i,j=1, \dots, l}}
\prod\limits_{1\leq i<j\leq l}\left(\frac{p_i-p_j}{q_i-q_j}\right)^2
\prod\limits_{i=1}^l \frac{ \rho^{\Pi}(p_i)}{\rho^{\Pi}(q_i)} {\Psi}^{\Pi}_{p_i,q_i}(X);
$$

2. For $\Prob_{\Pi}$-almost any $X\in
\Conf(U)$, the measure $\mathbb{P}(\cdot | X; U\setminus I)$ has the form
\begin{equation}
Z(I,X)^{-1} \prod\limits_{1\leq i<j\leq \#_I(X)} (t_i-t_j)^2 \prod \limits_{i=1}^{\#_I(X)}\rho^{\Pi}_{I,X}(t_i),
\end{equation}
where $Z(I,X)$ is the normalization constant
and the function $\rho^{\Pi}_{I,X}$ satisfies, for any $p,q\in I$, the relation
\begin{equation}
\frac{\rho^{\Pi}_{I,X}(p)}{\rho^{\Pi}_{I,X}(q)}=\frac{\rho^{\Pi}(p)}{\rho^{\Pi}(q)}\prod\limits_{x\in X\setminus I} \left(\frac{x-p}{x-q}\right)^2.
\end{equation}

\end{theorem}
{\bf Remark.} The order of the points in Claim 1 is  immaterial:  for any permutation $\pi$ on $l$ symbols, by definition, we  have 
$
\prod\limits_{i=1}^l {\Psi}^{\Pi}_{p_i,q_i}=
\prod\limits_{i=1}^l {\Psi}^{\Pi}_{p_i,q_{\pi(i)}}.
$

Let $U=(0, +\infty)$, take $s>-1$ and consider the  Bessel kernel
\begin{equation}\label{bessel}
J_s(x,y)=\frac{\sqrt{x}J_{s+1}(\sqrt{x})J_s(\sqrt{y})-\sqrt{y}J_{s+1}(\sqrt{y})J_s(\sqrt{x})}{2(x-y)}
\end{equation}
(see, e.g., page 295 in Tracy and Widom \cite{TracyWidom}). The kernel $J_s$ induces on $L_2((0, +\infty), {\mathrm{Leb}})$ the operator of orthogonal projection onto the subspace of functions whose
 Hankel transform is supported in $[0,1]$ (see \cite{TracyWidom}).

\begin{proposition}\label{rho-bessel}
For any $s>-1$, we have $\rho^{J_s}(t)=t^s$.
\end{proposition}

\subsection{The Hilbert-Schmidt Case.}
We now impose a weaker 
\begin{assumption}\label{as-key}
We have 
$
\displaystyle \int\limits_U \frac{\Pi(x,x)dx}{1+x^2}<+\infty.
$
\end{assumption}
 It follows that the operator
$
(x+i)^{-1}\Pi
$
is Hilbert-Schmidt.
The sine-kernel, for example, satisfies Assumption \ref{as-key} but not Assumption \ref{as-trace}.

Let $\la(x)$ be a continuous function on ${\mathbb R}$ satisfying
\begin{equation}\label{la-def}
\sup\limits_{x\in {\mathbb R}} \left| x^2\la(x)-x\right|<+\infty.
\end{equation}
For example, one can take $ \la(x)=(x+i)^{-1}$
or $
\la(x)=\displaystyle \frac{x}{x^2+1}$.

 We start by formulating an auxiliary
\begin{proposition}\label{la-reg}
\begin{enumerate}
\item 
For $p,q\in U$, the limit
\begin{equation}\label{mult-la}
{\Psi}^{\Pi, \la}_{p,q}(X)=
\lim\limits_{R\to\infty} \exp\left(2(p-q)\int\limits_{[-R, R]\cap U}  \Pi(x,x)\la(x)dx
\right) \prod\limits_{x\in X: |x|\leq R}
\left(\frac{x-p}{x-q}\right)^2
\end{equation}
exists in $L_1( \Conf(U), \Prob_{\Pi^q})$.
Furthermore,  for any compact subset $K\subset U$, there exists a subsequence $R_n\to\infty$, along which the almost sure convergence in (\ref{mult-la}) takes place for all $p,q\in K$.

\item There exists a positive function $\rho^{\Pi, \la}: U\to {\mathbb R}$ such that 
\begin{equation}\label{palm-rho-la}
\int \limits_{\Conf(U)} {\Psi}^{\Pi, \la}_{p,q}(X)d\Prob_{\Pi^q}(X)=\frac{ \rho^{\Pi, \la}(q)}{\rho^{\Pi, \la}(p)}\frac{\Pi(p,p)}{\Pi(q,q)}.
\end{equation}
\end{enumerate} 
\end{proposition}

If a configuration $X$ is represented in the form $X=\{t_1, \dots, t_l\}\cup Y$, where
$Y\in \Conf(U)$,
then, by definition,  we have
$$
{\Psi}^{\Pi, \la}_{p,q}(X)=\prod\limits_{i=1}^l\left(\frac{t_i-p}{t_i-q}\right)^2{\Psi}^{\Pi, \la}_{p,q}(Y).
$$

We are now ready to formulate the analogue of Theorem \ref{gibbs-trace}.
\begin{theorem}\label{gibbs-hs}
Let $U\subset {\mathbb R}$ be an open set. 
Let $\Pi(x,y)$, $x,y\in U$,  be a smooth kernel that 
 induces
 an operator of orthogonal projection acting in $L_2(U, {\mathrm {Leb}})$, satisfying Assumptions
\ref{rn-xpxq}, \ref{as-key} and such that the determinantal point process $\Prob_{\Pi}$ is rigid. 
Let $I\subset U$ be a compact interval.
Let $\la(x)$ be a continuous function on ${\mathbb R}$ satisfying
(\ref{la-def}).
Then 

1.  For almost any $2l$ distinct points $p_1, \dots, p_l, q_1, \dots, q_l\in U$, we have the  $\Prob^{q_1, \dots, q_l}$-almost sure equality
 $$
\frac {d\Prob^{p_1, \dots, p_l}}{d\Prob^{q_1, \dots, q_l}}(X)=
\frac {\det \Pi(q_i, q_j)|_{i,j=1, \dots, l}}{\det \Pi(p_i, p_j)|_{i,j=1, \dots, l}}
\prod\limits_{1\leq i<j\leq l}\left(\frac{p_i-p_j}{q_i-q_j}\right)^2
\prod\limits_{i=1}^l \frac{ \rho^{\Pi, \la}(p_i)}{\rho^{\Pi, \la}(q_i)} {\Psi}^{\Pi, \la}_{p_i,q_i}(X).
$$

2. For $\Prob_{\Pi}$-almost every $X\in
\Conf(U)$, the  measure $\mathbb{P}_{\Pi}(\cdot | X; U\setminus I)$ has the form
\begin{equation}\label{expform-cond-la}
Z(I,X, \la)^{-1} \prod\limits_{1\leq i<j\leq \#_I(X)} (t_i-t_j)^2 \prod \limits_{i=1}^{\#_I(X)}\rho^{\Pi, \la}_{I,X}(t_i),
\end{equation}
where $Z(I,X, \la)$ is the normalization constant
and the function $\rho^{\Pi, \la}_{I,X}$ satisfies, for any $p,q\in I$, the relation
\begin{equation}\label{rhopix-la}
\frac{\rho^{\Pi, \la}_{I,X}(p)}{\rho^{\Pi, \la}_{I,X}(q)}=\frac{\rho^{\Pi, \la}(p)}{\rho^{\Pi, \la}(q)}
{\Psi}^{\Pi, \la}_{p,q}(X|_{{\mathbb R}\setminus I}).
\end{equation}

\end{theorem}

{\bf Remark.}  1. The order of the points in Claim 1 is of course again  immaterial: see the Remark to Theorem \ref{gibbs-trace}.

2. Different choices of the function $\la$ result in the multiplication of ${\Psi}^{\Pi, \la}_{p,q}(X)$ by a constant.
More precisely, given continuous functions $\la_1$  and $\la_2$ satisfying (\ref{la-def}), the integral
$$
\beta_{\Pi}(\la_1, \la_2)=\int\limits_U (\la_1(x)-\la_2(x))\Pi(x,x) dx
$$
converges absolutely by Assumption \ref{as-key}. From the definitions we now have
$
 { \Psi}^{\Pi, \la_1}_{p,q}(X)= { \Psi}^{\Pi, \la_2}_{p,q}(X)\exp(2(p-q)\beta_{\Pi}(\la_1, \la_2)),
$
and, consequently, we have
$
\displaystyle \frac{\rho^{\Pi, \la_1}(p)}{\rho^{\Pi, \la_1}(q)}=
\displaystyle \frac{\rho^{\Pi, \la_2}(p)}{\rho^{\Pi, \la_2}(q)}\exp(2(q-p)\beta_{\Pi}(\la_1, \la_2)).
$
The expression (\ref{expform-cond-la}) does not, of course, depend on the specific choice of $\la$.

3.  Claim 2 of Theorem \ref{gibbs-hs} 
implies that for $\Prob_{\Pi}$-almost every $X\in \Conf(U)$ and any Borel automorphism $F$ of $U$ preserving the Lebesgue measure class and  acting by the identity 
in the complement of a compact subset $V\subset U$, setting $X\cap V=\{p_1, \dots, p_l\}$ and keeping the same symbol $F$ for the natural induced action of $F$ on the space of configurations, we have 
\begin{multline}\label{qi-hs}
\frac{d\Prob_{\Pi}\circ F}{d\Prob}\left(X\right)=\\=\prod\limits_{1\leq i<j\leq l}\left(\frac{F(p_i)-F(p_j)}{p_i-p_j}\right)^2 
\displaystyle \prod\limits_{i=1}^l \frac{ \rho^{\Pi, \la}(F(p_i))}{\rho^{\Pi, \la}(p_i)} \displaystyle  \frac{d\Leb\circ F}{d\Leb}(p_i){\Psi}^{\Pi, \la}_{F(p_i),p_i}(X|_{U\setminus V}).
\end{multline}

 Let $
{\mathscr S}(x,y)=\displaystyle \frac{\sin \pi(x-y)}{\pi(x-y)}
$
be the sine-kernel. For $\la_0(x)=x/(x^2+1)$ (any  odd function satisfying \eqref{la-def} would work), we have
\begin{equation}\label{psi-s}
\Psi^{{\mathscr S}, \la_0}_{p,q}(X)=\lim\limits_{R\to\infty} \prod\limits_{|x|\leq R} \left(\frac{x-p}{x-q}\right)^2.
\end{equation}
Convergence in (\ref{psi-s}) is in $L_1$ and almost sure along a subsequence, for instance, $R_n=n^4$.
Approximating the sine-kernel by Christoffel-Darboux kernels of Hermite polynomials in the usual way,
we obtain $\rho^{{\mathscr S}, \la_0}=1$. Theorem \ref{gibbs-hs} now yields
\begin{corollary}\label{rho-sine}
Let $I$ be a compact interval on ${\mathbb R}$.
For $\Prob_{\mathscr S}$-almost any configuration $X\in \Conf({\mathbb R})$, the conditional measure $\mathbb{P}_{{\mathscr S}}(\cdot | X; {\mathbb R}\setminus I)$ has the form
\begin{equation}
Z(I,X)^{-1} \prod\limits_{1\leq i<j\leq \#_I(X)} (t_i-t_j)^2 \prod \limits_{i=1}^{\#_I(X)}\rho^{\mathscr S}_{I,X}(t_i),
\end{equation}
where $Z(I,X)$ is the normalization constant
and the function $\rho^{{\mathscr S}}_{I,X}$ satisfies, for any $p,q\in I$, the relation
\begin{equation}\label{rhopix-sin}
\frac{\rho^{\mathscr S}_{I,X}(p)}{\rho^{{\mathscr S}}_{I,X}(q)}=\lim\limits_{R\to\infty}  \prod\limits_{x\in X\setminus I: |x|\leq R}
\left(\frac{x-p}{x-q}\right)^2.
\end{equation}

\end{corollary}

\section{Multiplicative functionals and Palm measures.}
\subsection{
Proof of Proposition \ref{la-reg}. }

Let $D_2\Pi$ stand for the Hessian of the kernel $\Pi$. The symbol $||\cdot||$ stands for the Euclidean norm of a vector or  a matrix.
\begin{lemma}\label{subst-est}
For any $\varepsilon>0$  and compact subset $K\subset U$,  there exists a positive constant $C(\varepsilon, K)$
such that for any $p,q\in K$ we have
$$
\sup\limits_{R\in {\mathbb R}} \left| \int\limits_{[-R, R]\cap U} \left(\left(\left(\frac{x-p}{x-q}\right)^2-1\right) \Pi^q(x,x)
+2(p-q) \Pi(x,x)\la(x)\right)dx
\right|\leq 
$$
$$
\leq C(\varepsilon, K)\left(1+\max\limits_{|x-q|\leq \varepsilon, |y-q|\leq \varepsilon} (||D_2\Pi||+||\Pi||)
 +\int\limits_U\frac{\Pi(x,x)dx}{1+x^2}\right).
$$
\end{lemma}
The proof of Lemma \ref{subst-est} is routine.
We represent the integral from  $-R$ to $R$ as a sum of two: first, the integral from $q-\varepsilon$ to $q+\varepsilon$, and, second, the integral over the remaining arcs. The first integral is estimated above by
$$C(\varepsilon, K)\max\limits_{|x-q|\leq \varepsilon, |y-q|\leq \varepsilon} (||D_2\Pi||+||\Pi||),$$
 the second, in view of (\ref{la-def}), by $C(\varepsilon, K)\displaystyle \int\limits_U\frac{\Pi(x,x)dx}{1+x^2}.$ The lemma is proved.

The result of \cite{Buf-gibbs} on the regularization of multiplicative functionals can be reformulated as follows:
\begin{lemma}\label{reg-mult}
For $p,q\in U$, the limit
$$
\lim\limits_{R\to\infty} \exp\left(-\int_{[-R,R]\cap U} \left(\left(\frac{x-p}{x-q}\right)^2-1\right) \Pi^q(x,x)dx
\right) \prod\limits_{x\in X: |x|\leq R}
\left(\frac{x-p}{x-q}\right)^2
$$
exists in $L_1( \Conf(U), \Prob_{\Pi^q})$.
Furthermore,  for any compact subset $K$ of $U$, there exists a subsequence $R_n\to\infty$, along which the almost sure convergence takes place for all $p,q\in K$.
\end{lemma}

Lemmata \ref{subst-est} and \ref{reg-mult} imply Proposition \ref{la-reg}.

\subsection{The function $\rho^{\Pi, \la}$.}

By Proposition \ref{la-reg}, we have
$
  { \Psi}^{\Pi, \la}_{p,q}(X)\in L_1(\Conf(U), \Prob_{\Pi^q}).
$
Assumption \ref{rn-xpxq} implies the relation
$$
L(p)=\frac{x-p}{x-q}L(q).
$$
By Corollary 4.12 in \cite{Buf-gibbs}, for any $p,q\in U$ there exists a positive constant $C_{\la}(p,q)$ such that for $\Prob^q$-almost every $X\in \Conf(U)$ we have
\begin{equation}\label{rn-der}
\frac{d\Prob_{\Pi}^{p}}{d\Prob_{\Pi}^{q}}(X)=C_{\la}(p,q)  \Psi_{pq}^{\Pi, \la}(X).
\end{equation}
For $p,q,r\in U$, we have
$
 \Psi_{pq}^{\Pi, \la} \Psi_{qr}^{\Pi, \la}= \Psi_{pr}^{\Pi, \la}
$
and 
$
C_{\la}(p,q)C_{\la}(q,r)=C_{\la}(p,r).
$

We now introduce a positive function $\rho^{\Pi, \la}$ on $U$ by setting
$$
C_{\la}(p,q)=\frac{\rho^{\Pi, \la}(p)\Pi(q,q)}{\rho^{\Pi, \la}(q)\Pi(p,p)},
$$
and (\ref{palm-rho-la}) is established. The function $\rho^{\Pi, \la}$ is of course defined up to a multiplicative constant.

In the case when the kernel $\Pi$ satisfies the stronger assumption (\ref{as-trace}), we can simply take $\la=0$ (even though $\la=0$ does not satisfy
(\ref{la-def})):  the operator $(x-q)^{-1}\Pi^q$ belongs to the trace class (since so does  $(x+i)^{-1}\Pi$), and we arrive at (\ref{palm-rho-pi}).

\subsection{Relation between Radon-Nikodym derivatives of Palm measures of different orders.}

As before, let $\Prob$ be a point process on a locally compact metric space $E$ endowed with 
a sigma-finite measure $\mu$ without atoms. As usual, we assume that for any $l$ the process $\Prob$ admits the $l$-th
 correlation measure of the form $\rho_l(p_1, \dots, p_l)d\mu(p_1)\dots d\mu(p_l)$.

\begin{proposition}\label{palm-dif}
Assume that for any natural number $l$ and $\mu^{\otimes l}$-almost  any two $l$-tuples $(p_1, \dots, p_l)$, $(q_1, \dots, q_l)$ 
of distinct points in $E$, the reduced Palm measures   $\Prob^{p_1, \dots, p_l}$ and $\Prob^{q_1, \dots, q_l}$ are equivalent.
Then for $\mu^{\otimes 2l}$-almost  any $2l$-tuple $(p_1, \dots, p_l, q_1, \dots, q_l)$ 
of distinct points in $E$ we have
$$
	\frac{\rho_l(p_1,...,p_l) d\mathbb{P}^{p_1,...,p_l}}{\rho_l(q_1,...,q_l) d\mathbb{P}^{q_1,...,q_l}}(X) 
= \prod\limits_{i=1}^l \frac{\rho_1(p_i)}{\rho_1(q_i)} \cdot \frac{d\mathbb{P}^{p_i}}{d\mathbb{P}^{q_i}}(X \cup q_1 \cup...\cup q_{i-1} \cup p_{i+1} \cup ... \cup p_l).
$$
\end{proposition}

\noindent {\bf Proof.} For $\mu$-almost any distinct $p, q, r_1, \dots, r_m\in E$,  we   clearly have
$$
	\frac{\rho_{m+1}(p,r_1,...,r_m)}{\rho_{m+1}(q,r_1,...,r_m)} \frac{d\mathbb{P}^{p,r_1,...,r_m}}{d\mathbb{P}^{q,r_1,...,r_m}}(X) = \frac{\rho_1(p)}{\rho_1(q)}\cdot \frac{d\mathbb{P}^p}{d\mathbb{P}^q}(X\cup r_1 \cup ... \cup r_m).
$$
The proposition is now proved by induction. For $l=2$ and $\mu$-almost any $p_1, p_2, q_1, q_2$, 
we have
$$
	\frac{\rho_2(p_1,p_2)}{\rho_2(q_1,q_2)}\frac{d\mathbb{P}^{p_1,p_2}}{d\mathbb{P}^{q_1,q_2}}(X) 
= \frac{\rho_2(p_1,p_2)}{\rho_2(q_1,p_2)}\frac{d\mathbb{P}^{p_1,p_2}}{d\mathbb{P}^{q_1,p_2}}(X)\cdot
\frac{\rho_2(q_1,p_2)}{\rho_2(q_1,q_2)}\frac{d\mathbb{P}^{q_1,p_2}}{d\mathbb{P}^{q_1,q_2}}(X) =
$$
$$
= \frac{\rho_1(p_1)}{\rho_1(q_1)} \frac{d\mathbb{P}^{p_1}}{d\mathbb{P}^{q_1}}(X\cup p_2) \cdot
\frac{\rho_1(p_2)}{\rho_1(q_2)} \frac{d\mathbb{P}^{p_2}}{d\mathbb{P}^{q_2}}(X\cup q_1).
$$

For  the induction step, we write 
$$
\frac{\rho_l(p_1,...,p_l) d\mathbb{P}^{p_1,...,p_l}}{\rho_l(q_1,...,q_{l-1},p_l) d\mathbb{P}^{q_1,...,q_{l-1},p_l}}(X)=\frac{\rho_l(p_1,...,p_{l-1}) d\mathbb{P}^{p_1,...,p_{l-1}}}{\rho_l(q_1,...,q_{l-1}) d\mathbb{P}^{q_1,...,q_{l-1}}}(X\cup p_l),
$$
whence, using the induction hypothesis, we conclude
$$
	\frac{\rho_l(p_1,...,p_l) d\mathbb{P}^{p_1,...,p_l}}{\rho_l(q_1,...,q_l) d\mathbb{P}^{q_1,...,q_l}}(X) = 
$$
$$
= 
\frac{\rho_l(p_1,...,p_l) d\mathbb{P}^{p_1,...,p_l}}{\rho_l(q_1,...,q_{l-1},p_l) d\mathbb{P}^{q_1,...,q_{l-1},p_l}}(X) \cdot \frac{\rho_l(q_1,...,q_{l-1},p_l) d\mathbb{P}^{q_1,...,q_{l-1},p_l}}{\rho_l(q_1,...,q_l) d\mathbb{P}^{q_1,...,q_l}}(X) =\
$$
$$
	=  \prod_{i=1}^{l-1}\frac{\rho_1(p_i)}{\rho_1(q_i)} \cdot \frac{d\mathbb{P}^{p_i}}{d\mathbb{P}^{q_i}}(X \cup q_1 \cup...\cup q_{i-1} \cup p_{i+1} \cup ... \cup p_{l}) 
\times \frac{\rho_1(p_l)}{\rho_1(q_l)} \cdot \frac{d\mathbb{P}^{p_l}}{d\mathbb{P}^{q_l}}(X \cup q_1 \cup ... \cup q_{l-1}).
$$
The proposition is proved completely.

\begin{corollary}\label{phi-palm}
Let $\Prob$ be a point process satisfying all assumptions of Proposition \ref{palm-dif}. If there exists a positive Borel function $\Psi: E \times E \times \Conf(E) \to \mathbb{R}_+$ and a positive Borel function $\Phi: \Conf_2(E) \to \mathbb{R}_+$ such that
\begin{enumerate}
	\item for $\mu$-almost any $p,q \in E,$  for $\mathbb{P}^q$-almost any $X\in \Conf(E)$, any $l \in \mathbb{N}$ and any distinct particles $r_1,...,r_l \in X$, we have
\begin{equation}\label{psi-phi}
	\Psi(p,q,X) = \frac{\Phi(p,r_1)}{\Phi(q,r_1)} \cdot \frac{\Phi(p,r_2)}{\Phi(q,r_2)} \cdot ... \cdot \frac{\Phi(p,r_l)}{\Phi(q,r_l)} \times \Psi(p,q,X\backslash\{r_1,...,r_l\});
\end{equation}
	\item for any $p, q, r \in E$ and $\mathbb{P}^q$-almost any $X\in \Conf(E)$ we have
$$
	\Psi(p,q,X) \cdot \Psi(q,r,X) = \Psi(p,r,X);
$$
	\item for $\mu$-almost any $p,q \in E$ and $\mathbb{P}^q$-almost any $X \in \Conf(E)$ we have
$$
	\frac{\rho_1(p)d\mathbb{P}^p}{\rho_1(q)d\mathbb{P}^q}(X) = \Psi(p,q,X),
$$
\end{enumerate} 	
then, for $\mu^{\otimes l}$-almost any $(p_1,...,p_l) \in \Conf_l(E),$ $(q_1,...,q_l)\in \Conf_l(E)$ and $\mathbb{P}^{q_1,...,q_l}$-almost 
any $X\in \Conf(E)$, we have 
$$
	\frac{\rho_l(p_1,...,p_l) d\mathbb{P}^{p_1,...,p_l}}{\rho_l(q_1,...,q_l)\mathbb{P}^{q_1,...,q_l}}(X) =
	\prod\limits_{1\le i<j \le l} \frac{\Phi(p_i,p_j)}{\Phi(q_i,q_j)} \cdot \prod\limits_{i=1}^l \Psi(p_i,q_i,X).
$$
\end{corollary}

Proposition \ref{palm-rho-la} , together with Proposition \ref{palm-dif} and Corollary \ref{phi-palm}, applied to our 
functional $\Psi_{p,q}^{\Pi, \la}$ satisfying (\ref{psi-phi}) with $\Phi(p,q)=|p-q|^2$, 
directly 
 implies the first claim of Theorems \ref{gibbs-trace}, \ref{gibbs-hs}.
We proceed to proving the second one. 

\subsection{Conditional Campbell measures.}
The following Proposition \ref{palm-cond}  will not be used in the proof and is included  to clarify the context. 

	Let $\mathbb{P}$ be a point process with locally finite intensity (in other words, admitting the first correlation measure) on $E$. Write $\xi\Prob$ for the first correlation measure of $\Prob$. Let $C \subset E$ be a Borel subset. 
	Let $\overline{\mathbb{P}_C}$ stand for the image of $\mathbb{P}$ under the natural projection map $\pi_C: \mathrm{Conf}(E) \to \mathrm{Conf}(C)$.

\begin{proposition}\label{palm-cond}
Assume that for $\mathbb{P}$-almost every $X$ the intensity $\xi\mathbb{P}_{(\cdot|X;C)}$ of the conditional process is absolutely continuous with respect to $\xi\mathbb{P}$. Then 
\begin{enumerate}
\item  for $\xi\mathbb{P}$-almost every $q\in E$ and $\overline{\mathbb{P}_C}$-almost every $Y \in \mathrm{Conf}(C)$ we have  $
	(\mathbb{P}^q)_{(\cdot|Y;C)}= (\mathbb{P}_{(\cdot|Y;C)})^q;
$  
\item for $\xi\mathbb{P}$-almost every $q \in E$ we have
\begin{equation}\label{palm-dec}
	\mathbb{P}^q= \int\limits_{\mathrm{Conf}(C)} \mathbb{P}_{(\cdot|Y;C)}^q \cdot \frac{d\xi\mathbb{P}_{(\cdot|Y;C)}}{d\xi\mathbb{P}}(q) \cdot d\overline{\mathbb{P}_C}(Y).
\end{equation}
\end{enumerate}
\end{proposition}

{\bf Proof.} Recall that the Campbell measure $\mathscr{C}_\mathbb{P}$ of the point process $\mathbb{P}$ is defined, for a compact  subset $B\subset E$ and a Borel subset $Z \subset \mathrm{Conf}(E)$, by the formula
$$
	\mathscr{C}_\mathbb{P} (B \times Z) = \int\limits_Z \#_B(X) \cdot d\mathbb{P}(X).
$$
By definition, we have
$
	\mathscr{C}_\mathbb{P} = \displaystyle \int\limits_{\mathrm{Conf}(C)} \mathscr{C}_{\mathbb{P}_{(\cdot|Y;C)}} d\overline{\mathbb{P}_C}(Y).
$
Let $\check \Prob^q$ stand for the {\it non-reduced} 
Palm measure of $\Prob$ at the point $q$.
We have
$
	\mathscr{C}_\mathbb{P} = \displaystyle \int\limits_E  {\check \Prob}^q \, d\xi\mathbb{P}(q)
$
and, similarly,
$
	\mathbb{P}_{(\cdot|Y;C)} = \displaystyle \int\limits_E {\check \Prob}_{(\cdot|Y;C)}^q \, d\xi\mathbb{P}_{(\cdot|Y;C)}(q).
$
Removing the point at $q$ and passing to reduced Palm measures,
we arrive at  \eqref{palm-dec}.

\begin{corollary}	\label{palm-cond-form}
Let $\Prob$ be a point process on $E$ such that for $\mathbb{P}$-almost every $X$ the intensity $\xi\mathbb{P}_{(\cdot|X;C)}$ of the conditional process is absolutely continuous with respect to $\xi\mathbb{P}$ and
for $\xi\mathbb{P}$-almost any $p,q \in E$ the reduced Palm measures $\Prob^p$ and $\Prob^q$ are equivalent.
Then for $\xi\mathbb{P}$-almost any $p,q \in E$ and $\mathbb{P}$-almost any $X\in \mathrm{Conf}(E)$ we have
\begin{equation}\label{palm-cond-rel}
	\frac{d\mathbb{P}^p}{d\mathbb{P}^q}(X) = \frac{\frac{d\xi\mathbb{P}_{(\cdot|X;C)}}{d\xi\mathbb{P}}(p)}{\frac{d\xi\mathbb{P}_{(\cdot|X;C)}}{d\xi\mathbb{P}}(q)} \cdot \frac{d\mathbb{P}_{(\cdot|X;C)}^p}{d\mathbb{P}_{(\cdot|X;C)}^q}(X|_{E\backslash C} ).
\end{equation}
\end{corollary}

Corollary \ref{palm-cond-form} is insufficient for our purposes: we need relation (\ref{palm-cond-rel}) to hold  on a fixed subset of  $\Conf(E)$ of full measure and for $\xi\mathbb{P}$-almost any $p,q \in E$.
To check this, we  use the quasi-invariance of our point processes under the group of compactly supported piecewise isometries of $E$.

\section{Palm measures and conditional measures.}
\subsection{Characterization of conditional measures.}
In this subsection, a general result is formulated describing conditional measures 
of point processes in terms of Radon-Nikodym derivatives of Palm measures of the same order.

Let $E$ be an open subset of ${\mathbb R}^d$, endowed with the Lebesgue measure $dv=dv_1\dots dv_d$. Let $\Prob$ be  a point process on $E$ satisfying   the following. 

\begin{assumption}\label{cor}
The point process $\Prob$ admits correlation measures of all orders.
For any $l>0$, the $l$-th correlation measure of $\Prob$ has the form
$$
\rho_l(p_1, \dots, p_l) dp_1\dots dp_l,
$$
where $\rho_l$ is a symmetric continuous function on $E^l$.
\end{assumption}

Recall that {\it the tail sigma-algebra} on $\Conf(E)$ is the intersection of all sigma-algebras 
$\mathcal{F}_{E\backslash B}$ over all compact $B\subset E$. 

\begin{assumption}\label{peq}
There exists a Borel subset ${\mathscr W}\subset \Conf(E)$, belonging to the tail sigma-algebra of $\Conf(E)$,
and, for any  $l>0$, a Borel measurable function $\Psi(p_1, \dots, p_l; q_1, \dots, q_l; X)$, defined
 for $X\in {\mathscr W}$ and  any two distinct $l$-tuples of points not containing particles of the configuration $X$, such that the following holds:
\begin{enumerate}
\item $\Prob({\mathscr W})=1$;
\item for fixed $X$, the function $\Psi(p_1, \dots, p_l; q_1, \dots, q_l; X)$ is continuous in $(p_1, \dots, p_l)\in \Conf_l(E\setminus X)$, $(q_1, \dots, q_l)\in \Conf_l(E\setminus X)$;
\item for fixed $X$ and  any three $l$-tuples $(p_1, \dots, p_l), (q_1, \dots, q_l), (r_1, \dots, r_l)$ in $\Conf_l(E\setminus X)$, we have
$$
\Psi(p_1, \dots, p_l; q_1, \dots, q_l; X)=\Psi(p_1, \dots, p_l; r_1, \dots, r_l; X)\Psi(r_1, \dots, r_l; q_1, \dots, q_l; X).
$$

\item for $\Prob$-almost any $Y\in {\mathscr W}$, any $l$ distinct particles $(p_1, \dots, p_l)\in  Y$ and
$\mu^{\otimes l}$-almost  any $l$-tuple
$(q_1, \dots, q_l)\in \Conf_l(E\setminus Y)$, we have
$$
\frac{d\Prob^{p_1, \dots, p_l}}{d\Prob^{q_1, \dots, q_l}}(Y)=\Psi(p_1, \dots, p_l; q_1, \dots, q_l; Y\setminus \{p_1, \dots, p_l\}).
$$
\end{enumerate}
\end{assumption}

\begin{proposition}\label{palm-cond-char}
Let $\Prob$ be a rigid point process on $E$ satisfying Assumptions \ref{cor}, \ref{peq}.
Let $I\subset E$ be a precompact open subset. Let $l\in {\mathbb N}$ be such that 
$$
\Prob(\{X: \#_I(X)=l\})>0.
$$
Then, for $\Prob$-almost every $X\in \Conf(E)$ such that $\#_I(X)=l$
and almost any distinct points $q_1, \dots, q_l \in E$,
the conditional measure $\Prob(\cdot |X, E\setminus I)$ has the form
\begin{equation}\label{form-cond}
Z_{q_1, \dots, q_l}^{-1}\Psi(p_1, \dots, p_l; q_1, \dots, q_l; X|_{E\setminus I}) \rho_l(p_1, \dots, p_l) dp_1\dots dp_l,
\end{equation}
 where
$Z_{q_1, \dots, q_l}$ is the normalization constant.
\end{proposition}

{\bf Remark.} The reference $l$-tuple $q_1, \dots, q_l \in E$ can be chosen arbitrarily; a different choice results in a change of the normalization constant.

\subsection{Quasi-invariance  under  piecewise isometries.}

We endow ${\mathbb R}^d$ with the norm $||v||=\max\limits_{i=1, \dots, d} |v_i|$ and the corresponding metric. 
The balls in this metric are cubes.
We take distinct points $p_1, \dots, p_l, q_1, \dots, q_l\in E$, take $\delta_1>0$, $\delta_2>0$, $\dots$, $\delta_l>0$ sufficiently small in such a way that the balls of radius $\delta_i$ centred at   $p_1, \dots, p_l, q_1, \dots, q_l$ do not intersect, and consider the piecewise isometry of $E$ that sends the closed ball of radius $\delta_i$
centred at $p_i$ to the corresponding ball centred at $q_i$, $i=1, \dots, l$, leaving the complement to the union of the closed balls fixed. The group generated by such piecewise isometries is denoted $\mathfrak G=\mathfrak G(E)$. For example,
 if $E={\mathbb R}$ , then the resulting group is the group of all interval exchange transformations on $\mathbb R$, while in higher dimension we arrive at the group of cube exchanges. 
The countable subgroup $\mathfrak G_0=\mathfrak G_0(E)$ generated by transformations of the above form such that the centres 
of all the balls have rational coordinates and the radii of the balls are rational.
For a subset $C\subset E$, let ${\mathfrak G}(C)$ and $\mathfrak G_0(C)$ be the subgroups of maps acting as the identity on  $E\setminus C$.  
For brevity, we write ${\bf p}=(p_1, \dots, p_l)$, 
$d{\bf p}=dp_1\dots dp_l$, $T{\bf p}=(Tp_1, \dots, Tp_l)$, etc.

	\begin{proposition} \label{qi-rn}
Let $I\subset {\mathbb R}^d$ be a bounded open set, let $l\in {\mathbb N}$. 
Let $F: \Conf_l(I) \to \mathbb{R}_+$ be a positive continuous function. Let $\mu$ be a Borel probability measure on $\Conf_l(I)$ such that the equality
	\begin{equation}\label{RN-der}
	\frac{d\mu \circ T}{d\mu}({\bf p}) = \frac{F(T{\bf p})}{F({\bf p})}.
	\end{equation}
 	holds $\mu$-almost surely for all $T \in \mathfrak G_0(I)$. Then (\ref{RN-der}) holds for all $T \in \mathfrak{G}(I)$ and 
$
d\mu({\bf p})=F({\bf p}) d{\bf p}.
$
	\end{proposition}

	{\bf Proof.} 
We first  show that $\mu$ assigns mass zero to boundaries of  balls:

\begin{lemma}\label{noatom}
For any $p\in I$ we have 
$
\mu(\{{\bf r}\in \Conf_l(I): p\in {\bf r}\})=0.
$
\end{lemma}

{\bf Remark.} The continuity of $F$ is essential, since any atomic measure with atoms of positive mass at all rational points 
in $\Conf_l(I)$ is quasi-invariant under $\mathfrak G_0(I)$.

{\bf Proof of Lemma \ref{noatom}.} 
First, we note that the measure $\mu$ cannot have atoms: if $\mu({\bf p})=\delta_0>0$, then, since 
the orbit of the configuration ${\bf p}$ under $\mathfrak G_0$ is dense in $\Conf_l(I)$ and
 (\ref{RN-der}) implies that there exists $\delta_1>0$ depending on $\delta_0$ and $F$
 such that the set $\{{\bf q}\in \Conf_l(I): \mu({\bf q})\geq \delta_1\}$ is infinite; but then the measure $\mu$ cannot be finite. Next, for any $i\leq d$ and any distinct points $p_1, \dots, p_i\in I$ we show 
$$
\mu(\{{\bf r}\in \Conf_l(I): p_1, \dots, p_i\in {\bf r}\})=0.
$$

We argue by induction on $i=d, d-1,\dots, 1$. The case $i=d$ is precisely the absence of atoms already established. 
For the induction step,  assume $\mu(\{{\bf r}: p_1, \dots, p_i\in {\bf r}\})>0$.  Then there exist points $q_1, \dots, q_i\in I$
and $\delta>0$, $\varepsilon>0$ and a ball $B(\varepsilon)$ of radius $\varepsilon$ in $\Conf_l(I)$ 
such that distances between distinct $q_k$ all exceed $2\varepsilon$   and we have
$\mu(\{{\bf r}: q_1, \dots, q_i\in {\bf r}\}\cap B(\varepsilon))>\delta$. Write $D=\{{\bf r}: q_1, \dots, q_i\in {\bf r}\}\cap B(\varepsilon)$. 
By  continuity of $F$, there exists $\delta_1>0$ such that the set of the ``shifts'' $TD$ of the set $D$ by elements $T\in 
\mathfrak G_0(I)$   satisfying $ \mu(TD)>\delta_1$ is infinite. The induction hypothesis implies $\mu(D\cap TD)=0$. 
It follows that the measure $\mu$ cannot be finite, and Lemma \ref{noatom} is proved completely.

We proceed with the proof of Proposition \ref{qi-rn}.
A ball of radius $r$ centred at a configuration ${\bf p}\in \Conf_l(I)$ will be called {\it proper} 
if the distances between the distinct $p_i$ are all less than $r/2$. 

Take two finite collections $B_1, \dots, B_k$, $B^{\prime}_1, \dots, B^{\prime}_k$ of disjoint isometric proper balls and let $T$ be a piecewise isometry interchanging $B_i$ and $B^{\prime}_i$, $i=1, \dots, k$. To establish Proposition \ref{qi-rn}, it suffices to establish 
(\ref{RN-der}) for piecewise isometries $T$ of this form.

	Take an exhausting sequence $B_{n,i}\subset B_i$, $B_{n,i}^{\prime}\subset B^{\prime}_i$ of isometric balls with rational centres and radii.   
Define $T_n\in \mathfrak G_0$ as the map that interchanges $B_{n,i}$ and $B_{n,i}^{\prime}$, $i=1, \dots, k$.
Lemma \ref{noatom} implies  that the sequence $\mu \circ T_n$ weakly converges to $\mu \circ T$ and also that the sequence $F(T_n{\bf p})\mu$ weakly converges to the limit $F(T{\bf p})\mu$ as $n\to \infty$.

	Take $\varepsilon > 0$. Set  $\Conf_{l, \varepsilon}(I)=\{{\bf p}\in \Conf_l(I): \min\limits_{i,j = 1,...,l} |p_i - p_j| \geq \varepsilon\} $. Let  $\varphi$ be a bounded continuous   
function  on $\Conf_l(I)$ supported on $\Conf_{l, \varepsilon}(I)$.
	The function ${\varphi({\bf p})}/{F({\bf p})}$ is  then bounded and continuous, and we have
$$
	\lim\limits_{n\to\infty} \int\limits_{\Conf_l(I)} \varphi({\bf p}) \cdot \frac{F(T_n{\bf p})}{F({\bf p})} d\mu({\bf p}) 
	= \int\limits_{\Conf_l(I)} \varphi({\bf p})\cdot \frac{F(T{\bf p})}{F({\bf p})} d\mu({\bf p}),
$$
whence  the sequence of probability measures $\displaystyle \frac{F(T_n{\bf p})}{F({\bf p})}\mu = \mu \circ T_n$ vaguely converges,  as $n\to\infty$, to the measure $\displaystyle \frac{F(T{\bf p})}{F({\bf p})}\mu$.
Since the sequence $\mu \circ T_n$ weakly converges to $\mu \circ T$, the equality \eqref{RN-der} is proved for all $T\in {\mathfrak G}$.

To conclude the proof of Proposition \ref{qi-rn}, consider the measure $\eta$ given by 
$d\eta({\bf p})=d\mu({\bf p})/F({\bf p})$. By continuity and positivity of $F$, for any $\varepsilon>0$, the measure $\eta$ is finite in restriction to $\Conf_{l, \varepsilon}(I)$. Since $\eta$ is ${\mathfrak G}$-invariant,  the measure $\eta$, in restriction to $\Conf_{l, \varepsilon}(I)$, coincides with the Lebesgue measure. Since $\varepsilon$ is arbitrary, Proposition \ref{qi-rn} is proved completely.

\subsection{ Completion of the proof of Proposition \ref{palm-cond-char}.}
	Let $S$ be a standard Borel space, let $\mu$ be a Borel probability measure on $S$. Let $\mathcal{F}$ be a $\sigma$-algebra of Borel subsets of $S$, let $\pi$ be the corresponding measurable partition. We let $\bar\mu$ be the quotient measure of $\mu$ under the partition $\pi$, and, for an element $\xi$ of the partition $\pi$, we let $\mu^\xi$ be the corresponding conditional measure. Finally, let $T$ be a Borel transformation of the space $S$ such that every set of $\mathcal F$ is $T$-invariant and the measure $\mu$ is $T$-quasi-invariant.
	The definitions directly imply 
\begin{proposition}\label{cond-qi}
	Let $F$ be a Borel function such that the equality
$$
	\frac{d\mu \circ T}{d\mu} = F
$$
 	holds $\mu$-almost surely. Then for $\bar\mu$-almost every element $\xi$ of the partition $\pi$ we have the $\mu^\xi$-almost sure equality
$$
	\frac{d\mu^\xi \circ T}{d\mu^\xi} = F.
$$
\end{proposition}

Proposition 2.9 in \cite{Buf-gibbs} claims that for 
a piecewise isometry $T\in \mathfrak G$ acting as the identity beyond 
a compact set $V$ and a configuration $X\in \Conf(E)$ such that $X\cap V=\{p_1, \dots, p_l\}$, 
we have, $\Prob$-almost surely, the equality
\begin{equation}\label{tqi}
\frac{d\Prob\circ T}{d\Prob}(X)=\frac {\rho_l(Tp_1, \dots Tp_l)}{\rho_l(p_1, \dots p_l)}\frac{d\Prob^{Tp_1, \dots, Tp_l}}{d\Prob^{p_1, \dots, p_l}}(X\setminus \{p_1, \dots, p_l\}).
\end{equation}
Let $I\subset E$ be  precompact and open.
By Proposition \ref{cond-qi}, for $\Prob$-almost any $X\in \Conf(E)$ and any $T\in 
\mathfrak G_0$, the measure $\Prob(\cdot|X, E\setminus I)$
satisfies the equality (\ref{tqi}) (in which one must, of course, substitute $\Prob(\cdot|X, E\setminus I)$ for $\Prob$). 
By Proposition \ref{qi-rn}, the same equality  holds for all $T\in \mathfrak G$ and the  measure $\Prob(\cdot|X, E\setminus I)$ has  the form 
(\ref{form-cond}). Proposition \ref{palm-cond-char} is proved completely.

\section{Continuity of the functions $\rho^{\Pi}$, $\rho^{\Pi, \la}$ and the proofs of Proposition \ref{rho-bessel}, Corollary \ref{rho-sine}.
}
\subsection{The trace class case}
Let $D_1\Pi$ stand for the Jacobi matrix of the kernel $\Pi$.
Our definitions immediately  imply the following important continuity property of the function $\rho^{\Pi}$. 
\begin{proposition}\label{tr-conv}
Let $\Pi_n$ be a sequence of smooth kernels, each inducing an operator of orthogonal projection in
$L_2(U, \Leb)$, each satisfying  Assumptions \ref{rn-xpxq} and  \ref{as-trace}. Assume that, as $n\to\infty$, we have
\begin{enumerate}
\item $\Pi_n \to \Pi, D_1\Pi_n\to D_1\Pi, D_2\Pi_n \to D_2\Pi
$
uniformly on compact subsets of $U\times U$ ;

\item $(|x|+1)^{-1}\Pi_n\to (|x|+1)^{-1}\Pi$ in the space 
of trace class operators acting in $L_2(U, \Leb)$.
\end{enumerate}
Then, for any any $p,q\in U$, we have
$
\lim\limits_{n\to\infty} \displaystyle \frac{ \rho^{\Pi_n}(p)}{\rho^{\Pi_n}(q)}=
\displaystyle \frac{\rho^{\Pi}(p)}{\rho^{\Pi}(q)}.
$
\end{proposition}

\subsection{The Bessel kernel: computation of the function $\rho^{J_s}$ .}

Let $s>-1$.  Let $P_n^{(s)}$ be the
standard Jacobi orthogonal polynomials corresponding to
the weight $(1-u)^s$. Let $\tilde K_n^{(s)} (u_1, u_2)$ the $n$-th Christoffel-Darboux kernel of the Jacobi orthogonal  polynomial ensemble.
Recall that the classical Heine-Mehler asymptotics for Jacobi orthogonal polynomials (see e.g.
Chapter 8 in Szeg{\"o} \cite{Szego}) implies that  for any $s>-1$, as $n\to\infty$, the kernel ${\tilde K}_n^{(s)}$
converges
to the kernel ${J}_{s}$ uniformly in the totality of variables
on compact subsets of
$(0, +\infty)\times (0, +\infty)$, indeed, 
 on arbitrary simply connected compact subsets of
$({\mathbb C}\setminus 0)\times ({\mathbb C}\setminus 0)$.
Our next aim is to justify the limit transition
\begin{equation}\label{tr-conv-jac}
\lim\limits_{n\to\infty} \frac{ \rho^{\Pi_n}(p)}{\rho^{\Pi_n}(q)}=\lim\limits_{n\to\infty} \frac{ (1-(1-p/2n^2))^s}{(1-(1-q/2n^2))^s}=\frac{p^s}{q^s}=
\frac{\rho^{J_s}(p)}{\rho^{J_s}(q)}.
\end{equation}
By  Proposition \ref{tr-conv}, it remains to   prove that 
$(1+x)^{-1}{\tilde K}_n^{(s)}\to (1+x)^{-1}{J}_{s}$ in the space of trace class 
operators acting in $L_2({\mathbb R}_+)$.  
For  $s>0$, this
trace class convergence directly follows  from  standard  inequalities for Jacobi polynomials, see as e.g. Theorem 7.3.2 in Szeg{\"o} \cite{Szego}.
To treat the case $s\in (-1,0]$, note that for any $s>-1$ we have the recurrence relations
\begin{multline} \label{rec-form-jac-s}
\tilde K_{n}^{(s)}(u_1, u_2)=
\frac{s+1}{2^{s+1}}P_{n-1}^{(s+1)}(u_1)(1-u_1)^{s/2} P_{n-1}^{(s+1)}(u_2)(1-u_2)^{s/2}+\\+ \tilde K_{n}^{(s+2)}(u_1, u_2)
\end{multline}
\begin{equation}\label{rec-bes}
J_s(x,y)=J_{s+2}(x,y)+\frac{s+1}{\sqrt{xy}}J_{s+1}(\sqrt{x})J_{s+1}(\sqrt{y}).
\end{equation}
Relations (\ref{rec-form-jac-s}), (\ref{rec-bes}) imply the  convergence $(1+x)^{-1}{\tilde K}_n^{(s)}\to (1+x)^{-1}{J}_{s}$ in trace class norm for any $s>-1$.
Proposition \ref{rho-bessel} is proved completely.

\subsection{The Hilbert-Schmidt Case}
Our definitions directly imply
\begin{proposition}\label{hs-conv}
Let $\Pi_n$ be a sequence of smooth kernels, each inducing an operator of orthogonal projection acting in
$L_2(U, \Leb)$, each satisfying  Assumptions \ref{rn-xpxq} and \ref{as-key}.
 If, as $n\to\infty$, we have
\begin{enumerate}
\item
 $\Pi_n \to \Pi, D_1\Pi_n\to D_1\Pi, D_2\Pi_n \to D_2\Pi
$
uniformly on compact subsets of $U\times U$ ;
\item $(x+i)^{-1}\Pi_n\to (x+i)^{-1}\Pi$ in the Hilbert-Schmidt norm,
\end{enumerate}
then, for any continuous $\la$ satisfying (\ref{la-def}) and any $p,q\in U$ we have
\begin{equation}
\lim\limits_{n\to\infty} \frac{ \rho^{\la}_{\Pi_n}(p)}{\rho^{\la}_{\Pi_n}(q)}=
\frac{\rho^{\la}_{\Pi}(p)}{\rho^{\la}_{\Pi}(q)}.
\end{equation}
\end{proposition}
\begin{proposition}\label{power-est-hs}
If $\Pi_n\to\Pi$  uniformly on compact subsets of $U$ and there exists $\alpha$, 
$0\leq \alpha <1/2$, such that
\begin{equation}\label{unif-power-hs}
\sup\limits_{n\in {\mathbb N}, x\in {\mathbb R} } \frac{\Pi_n(x,x)}{1+|x|^{\alpha}}<+\infty.
\end{equation}
Then $(x+i)^{-1}\Pi_n\to (x+i)^{-1}\Pi$
 in Hilbert-Schmidt norm.
\end{proposition}

Proof. Indeed, by Gr{\"u}mm's theorem (see e.g. Simon \cite{Simon}), it suffices
to check the relation
\begin{equation}\label{limtrace}
\lim\limits_{n\to\infty}\int\limits_{-\infty}^{\infty} \frac{\Pi_n(x,x)dx}{1+x^2}=\int\limits_{-\infty}^{\infty}\frac{\Pi(x,x)dx}{1+x^2}.
\end{equation}
For any $R_0>0$, the uniform convergence of our kernels on compact subsets
implies the convergence
$$
\lim\limits_{n\to\infty}\int\limits_{-R_0}^{R_0} \frac{\Pi_n(x,x)dx}{1+x^2}=\int\limits_{-R_0}^{R_0}\frac{\Pi(x,x)dx}{1+x^2}.
$$
Condition (\ref{unif-power-hs}), in turn, immediately implies, for any $\varepsilon>0$, the existence of
$R_0>0$ such that
\begin{equation}\label{unif-power-2-hs}
\sup\limits_{n\in {\mathbb N} } \int\limits_{|x|>R_0}\frac{\Pi_n(x,x)}{1+x^{2}}<+\infty,
\end{equation}
 convergence (\ref{limtrace}) follows, and Proposition \ref{power-est-hs} is proved.

\subsection{The sine-kernel}

Let $\la_0(x)=x(x^2+1)^{-1}$ so that (\ref{psi-s}) holds. 
Since
$$
{\mathrm {Var}}_{\Prob_{\mathscr S}} \left(\sum\limits_{x\in X, |x|\geq R} (\log|x-p|-\log|x-q|)\right) =O(R^{-2}).
$$
the Borel-Cantelli lemma implies convergence  in (\ref{psi-s}), for example, along the sequence  $R_n=n^4$.
Let
$
{\widetilde K}_n^{(H)}
$
be the Christoffel-Darboux kernel of the standard Hermite polynomials and set $$K^{(H)}_n(x,y)=\frac{\pi}{\sqrt{2n} }{\widetilde K}_n^{(H)}(\frac{x}{\sqrt{2n}}, \frac{y}{\sqrt{2n}}).$$
We have 
$
\lim\limits_{n\to\infty}K^{(H)}_n(x,y)= {\mathscr S}(x,y).
$
Convergence is uniform with all derivatives as long as $x,y$ range over compact subsets of the complex plane. The Plancherel-Rotach 
asymptotic for Hermite polynomials, see e.g. Theorem 8.22.9 in Szeg{\"o} \cite{Szego}, implies  (\ref{unif-power-hs}) for $\Pi_n=K^{(H)}_n$, and Proposition \ref{power-est-hs} implies the  Hilbert-Schmidt convergence $(x+i)^{-1}K^{(H)}_n\to (x+i)^{-1} {\mathscr S}$.

Since  $\la_0$ is odd and  $K^{(H)}_n(x,x)$ is even, similarly to (\ref{psi-s}),
we have
$$
\Psi^{K^{(H)}_n, \la_0}_{p,q}(x_1, \dots, x_n)=\prod\limits_{i=1}^n \left(\frac{x_i-p}{x_i-q}\right)^2.
$$

Since $\lim\limits_{n\to\infty}\exp(-p^2/2n+q^2/2n)=1,$
we conclude 
$
\rho^{{\mathscr S}, \la_0}(p)=1.
$

{\bf Remark. } The Airy kernel satisfies all assumptions of Theorem \ref{gibbs-hs}; the  explicit constants will be given in the sequel to this paper.

{\bf Acknowledgements.}
I am deeply grateful to Grigori Olshanski, Valery Oseledets, Sergei Pirogov, Yanqi Qiu and Dmitri Zubov for useful discussions.This project has received funding from the European Research Council (ERC) under the European Union's Horizon 2020 research and innovation programme (grant agreement No 647133 (ICHAOS))
and has also been funded by the Russian Academic Excellence Project `5-100'. 
I am deeply grateful to  Institut Henri Poincar{\'e}, where part of this work was carried out, for its warm hospitality.

\end{document}